\documentclass[12pt,letter]{article} 
\usepackage[T2A]{fontenc}
\usepackage[cp866]{inputenc}
\usepackage[russian,english]{babel}
\usepackage[tbtags]{amsmath}
\usepackage{amsfonts,amssymb,mathrsfs,amscd,comment}
\usepackage{amsthm,graphicx} 
\usepackage{hyperref}
\def\NoBlackBoxes{\overfullrule0pt}

\NoBlackBoxes
\theoremstyle{plain}

\newtheorem{conj}{Conjecture}
\theoremstyle{definition}

\newtheorem{remark}{Remark}

\theoremstyle{main}
\def\bad{\spaceskip=0.33emplus0.6emminus0.15em\immediate\write5{\string\bad}}

\theoremstyle{plain}
\newtheorem{proposition}{Proposition}
\let\myh\widehat
\let\myo\overline
\let\eps\varepsilon
\let\pfi\varphi
\def\HH{\mathscr H}
\def\RR{\mathbb R}
\def\CC{\mathbb C}
\def\QQ{\mathbb Q}
\def\NN{\mathbb N}
\def\ZZ{\mathbb Z}
\def\sZ{\mathscr Z}
\def\ext{\operatorname{ext}}
\def\supp{\operatorname{supp}}
\def\const{\operatorname{const}}
\def\<{\left\langle}
\def\>{\right\rangle}
\def\({\left(}
\def\){\right)}
\def\[{\left[}
\def\]{\right]}
\NoBlackBoxes

\begin{document}

\selectlanguage{english}


\author{Sergey~P.~Suetin}

\title{On one example of a~Nikishin system}

\markboth{S.~P.~Suetin}{On one example of a Nikishin system}

\maketitle


\begin{abstract}
The paper puts forward an example of a~Markov function $f=\const+\myh{\sigma}$ such that
the three functions $f,f^2$ and $f^3$ form a~Nikishin system. A~conjecture is proposed
that there exists a~Markov function $f$ such that, for each
$n\in\NN$, the system $f,f^2,\dots,f^n$ constitutes a~Nikishin system.

Bibliography:~20~titles.
\end{abstract}


\footnotetext{This research was carried out with the partial financial support of the Russian Foundation
for Basic Research (grant no.\ 18-01-00764).}

\section{Introduction and statement of the problem}\label{s1}
As distinct from Pad\'e polynomials, which are constructed from one
function\footnote{More precisely, here we speak about the power series
expansion defined
at some fixed
point $z_0$ of the Riemann sphere $\myo\CC$, for example, $z_0=\infty$.} $f$, a~construction of the
Hermite--Pad\'e polynomials corresponding, for example, to a two-\allowbreak dimensional 
multiindex,
requires at least two functions $f_1$ and $f_2$, which should be in a~sense independent.
Namely, in order that the definition of the Pad\'e polynomials be meaningful it is necessary that the
original function~$f$ should not be a~rational function. In other words,
it is necessary that the pair of functions $1,f$ should be independent over the field of rational functions $\CC(z)$.
Likewise, in order that the definition of the Hermite--Pad\'e polynomials for a~pair of functions $f_1,f_2$
be meaningful it is required that the three functions $1,f_1,f_2$ be independent over the field $\CC(z)$.
For the definition of Hermite--Pad\'e polynomials and general properties of these polynomials, see, above all,
\cite{Nut84} and~\cite{NiSo88}, and also~\cite{Van06}.

The family of functions involved in the construction of Hermite--Pad\'e poly\-no\-mi\-als 
is usually called a~{\it system}.
The two best known systems in the theory of Hermite--Pad\'e polynomials are an
{\it Angelesco system} and a~{\it Nikishin system};
for the formal definition of such systems and their properties, see, above all, \cite{NiSo88}, and also
\cite{GoRaSo97}, \cite{Van06}, \cite{FiLoLoSo10}, \cite{FiLo11}, \cite{Rak18}. More general
(the so-called ``mixed'') systems of functions were considered in \cite{Sor83}, \cite{ApLy10}, \cite{Rak11}.
Below we will briefly discuss the meaningfulness of these concepts in the case of a~pair of
functions $f_1$ and~$f_2$.

Given an arbitrary (positive Borel) measure $\sigma$ with support $\supp{\sigma}$
on the real line $\RR$, $\supp{\sigma}\Subset\RR$, we denote by
\begin{equation}
\myh{\sigma}(z):=\int\frac{d\sigma(x)}{z-x},\quad z\notin\supp{\sigma},
\label{01}
\end{equation}
the Cauchy transform of the measure $\sigma$.

For a pair of functions $f_1$ and $f_2$ of the form~\eqref{01} the property that this pair forms an {\it Angelesco system}
appears to be quite natural. Namely, in this case the functions $f_1$ and $f_2$ can be written as
\begin{equation}
f_1(z):=\myh{\sigma}_1(z),\quad
f_2(z):=\myh{\sigma}_2(z),
\label{1}
\end{equation}
where it is assumed that the supports of the measures $\sigma_1$ and $\sigma_2$ are disjoint,
$\supp{\sigma_1}\cap\supp{\sigma_2}=\varnothing$,
$\supp{\sigma_1},\supp{\sigma_2}\Subset\RR$.

If in \eqref{1} the supports of the measures $\sigma_1$ and $\sigma_2$ are equal,
$\supp{\sigma_1}=\supp{\sigma_2}=\Delta$, and if
\begin{equation}
d\sigma_2(x)=\myh{\sigma}_3(x)\,d\sigma_1(x),\quad x\in \Delta,
\label{1.1}
\end{equation}
where the third measure $\sigma_3$, $\supp{\sigma_3}\Subset\RR$, is such that
$\supp{\sigma}_3\cap \Delta=\varnothing$, then the pair of functions $f_1$ and~$f_2$ of the form \eqref{01} is said to
form a~{\it Nikishin system}.

At first glance the definition of an Angelesco system looks more natural than that of a~Nikishin system. For
example, an Angelesco system is formed by the pair of functions
\begin{equation}
f_1(z):=\frac1{\bigl[(z-a_1)(z-b_1)\bigr]^{1/2}},\quad
f_2(z):=\frac1{\bigl[(z-a_2)(z-b_2)\bigr]^{1/2}},
\label{3}
\end{equation}
where $a_1<b_1<a_2<b_2$ and a~branch of the root is chosen so that
$\bigl[(z-a_j)(z-b_j)\bigr]^{1/2}/z\to1$ as $z\to\infty$.

The purpose of the present note is to present an example of
a~Markov function $f=\myh{\sigma}+\const$ such that the pair of functions $f,f^2$ forms a~Nikishin system
under a~certain minimal extension of the original definition of such a~system
(see~\eqref{10}--\eqref{12}, and also Remark~\ref{rem2} below). As a~result, it
turns out that, from the point of view
of the general problem of efficient analytic continuation of a~multivalued analytic function
defined by a~power series (for more details, see~\cite{Sue18c}), the concept of a~Nikishin system
is by no means less meaningful
than that of an Angelesco system. It is worth noting that this fact is also manifested in
some
papers on Nikishin systems; see, for example, \cite{KoPaSuCh17},
\cite{ApBoYa17}, \cite{BaGeLo18}, \cite{LoMi18},
\cite{LoVa18}, \cite{Sue18},~\cite{Sue18b}
and the references given therein.

More precisely, we will give an example of a function $f$ of the form
\begin{equation}
f(z)=C+\myh{\sigma}(z)
\label{4}
\end{equation}
(cf.~\eqref{01}), where $C\neq0$ is some real constant, $\sigma$ is a~measure
supported on the interval
$[-1,1]$, $\supp\sigma=[-1,1]$, such that the pair of functions
\begin{equation}
f_1(z):=f(z),\quad f_2(z):=f^2(z)
\label{5}
\end{equation}
forms a~Nikishin system.
Furthermore, it will be shown that,
for the function $f(z)$ considered below (see~\eqref{8}) of the form \eqref{4},
the three functions $f,f^2,f^3$ also
form a~Nikishin system.

\begin{remark}\label{rem1}
One consequence of the presence of the term $C\neq0$ in representation~\eqref{4} is that
the function $f_2(z)=f^2(z)$ can no longer be written in the form \eqref{1}--\eqref{1.1}.
Nevertheless, somewhat more involved representations will be shown to hold. Namely, the following representations
are valid
\begin{equation}
f(z)-C=\myh{\sigma}(z),\quad
f^2(z)-Cf(z)=\myh{s}_1(z),\quad
f^3(z)-Cf^2(z)=\myh{s}_2(z),
\label{6.2}
\end{equation}
where
$\supp s_j=[-1,1]$, $j=1,2$; for more details, see \S~\ref{s2} and Remark~\ref{rem2} below.

The possibility of the existence of a Markov function~$f$ for which similar representations
would hold for an arbitrary power $f^n$ will be discussed below
(see Conjecture~\ref{con1}).
\end{remark}

\section{Definitions and statement of the main result}\label{s2}

Let $\Delta_1:=[-1,1]$,
\begin{equation}
\pfi(z):=z+(z^2-1)^{1/2}, \quad z\notin \Delta_1,
\label{6}
\end{equation}
be the function inverse to the Zhukovskii function. Recall that we have chosen and fixed
a~branch of the square root such that $(z^2-1)^{1/2}/z\to1$ as $z\to\infty$.
So, $|\pfi(z)|>1$ for $z\notin \Delta_1$.
Hence, for any complex number~$A$ such that $|A|>1$, the multivalued analytic function
\begin{equation}
f(z):=f(z;A,\alpha):=\(A-\frac1{\pfi(z)}\)^\alpha,
\quad\text{where}\quad\alpha\in\CC\setminus\ZZ,
\label{7}
\end{equation}
admits a~holomorphic (i.e., a~single-valued analytic)
branch in the domain $D_1:=\myo\CC\setminus \Delta_1$. However, in the domain
$\myo\CC\setminus\{-1,1\}$ the function
$f(z)$ is already a~multivalued analytic function, the set of branch points~$\Sigma$ of this
function consisting of three points: $\Sigma=\{\pm1,a\}$,
where $a=(A+1/A)/2$ and, hence, $|a|>1$. Note that $1/\pfi(z)=z-(z^2-1)^{1/2}$
in accordance with the above choice of the branch of the root in~\eqref{7}.

The class of multivalued analytic functions $\mathscr Z$ consisting of all functions obtained by
multiplication of a~finite number of functions of the form~\eqref{7}
\begin{equation}
f(z):=\prod_{j=1}^m\(A_j-\frac1{\pfi(z)}\)^{\alpha_j},
\label{7.2}
\end{equation}
where $|A_j|>1$, $\alpha_j\in\CC\setminus\ZZ$ for all $j=1,\dots,m$, and $\sum_{j=1}^m\alpha_j\in\ZZ$,
was introduced and studied in~\cite{Sue17} (see
also \cite{Sue17b},~\cite{Sue18c}).
In the present paper, we shall be concerned only with the case when in~\eqref{7.2}\enskip $m=2$,
$A_1,A_2$ are real, and $\alpha_1=\alpha_2=-1/2$. We shall also
assume that $1<A_1<A_2$. So, the functions to be considered are of the form
\begin{equation}
f(z):=\[\(A_1-\frac1{\pfi(z)}\)\(A_2-\frac1{\pfi(z)}\)\]^{-1/2},
\label{8}
\end{equation}
or, in other words, $f(z)=f_1(z)f_2(z)$, where
\begin{equation}
\begin{aligned}
f_1(z):&=f(z;A_1,-1/2):=\(A_1-\frac1{\pfi(z)}\)^{-1/2}
=\frac1{\(A_1-1/\pfi(z)\)^{1/2}},\\
f_2(z):&=f(z;A_2,-1/2):=\(A_2-\frac1{\pfi(z)}\)^{-1/2}
=\frac1{\(A_2-1/\pfi(z)\)^{1/2}},\\
\end{aligned}
\label{9}
\end{equation}
$z\in D_1$, $1<A_1<A_2$.

In what follows, $\sqrt{\cdot}$ denotes the positive square root
of a~nonnegative real number; i.e., $\sqrt{a^2}=|a|$ for $a\in\RR$.

The main result of the present paper is as follows.

\begin{proposition}\label{pro1}
Let $f(z)$ be the function defined by representation~\eqref{8}, where $1<A_1<A_2$, and let
 $a_j=(A_j+1/A_j)/2$, $j=1,2$. Then, for $z\in D$,
\begin{align}
f(z)&=\frac1{\sqrt{A_1A_2}}+\myh{\sigma}(z),
\label{10}\\
f^2(z)&=\frac1{A_1A_2}+\frac1{\sqrt{A_1A_2}}\myh{\sigma}(z)
+\myh{s}_1(z),
\label{11}\\
f^3(z)&=\frac1{\sqrt{(A_1A_2)^3}}
+\frac1{A_1A_2}\myh{\sigma}(z)+\frac1{\sqrt{A_1A_2}}\myh{s}_1(z)
+\myh{s}_2(z),
\label{12}
\end{align}
where $\sigma$ is the measure supported on the interval $[-1,1]$, the measures $s_1$ and $s_2$
are defined by the representations $s_1=\<\sigma,\sigma_2\>$
and $s_2=\<\sigma,\sigma_2,\sigma\>$. Moreover, $\supp{s_1}=\supp{s_2}=[-1,1]$,
$\supp{\sigma_2}=[a_1,a_2]\subset\RR\setminus[-1,1]$, and the measures $\sigma$ and~$\sigma_2$
have the following explicit representations
\begin{align}
d\sigma(x_1)&=
\frac{\sqrt{1-x^2_1}}{4\pi\sqrt{A_1A_2}\sqrt{(a_1-x_1)(a_2-x_1)}}
\[\frac{h_2(x_1)}{h_1(x_1)}+\frac{h_1(x_1)}{h_2(x_1)}\]\,dx_1,
\quad x_1\in[-1,1],
\label{13}\\
d\sigma_2(x_2)&=
\frac1{\pi}\frac{dx_2}{\sqrt{(\pfi(x_2)-A_1)(A_2-\pfi(x_2))}},\quad
x_2\in(a_1,a_2),
\label{14}
\end{align}
where
$$
h_j(x_1):=\(A_j-(x_1+i\sqrt{1-x^2_1})\)^{1/2}+\(A_j-(x_1-i\sqrt{1-x^2_1})\)^{1/2}>0
$$
for $x_1\in[-1,1]$, $j=1,2$.
\end{proposition}

Following \cite{GoRaSo97}, in Proposition \ref{pro1} we used the following notation
for the measure~$s_1$:
$d\<\sigma,\sigma_2\>(x_1):=\myh{\sigma}_2(x_1)\,d\sigma(x_1)$,
$x_1\in\Delta_1^\circ:=(-1,1)$, which is legitimate under our assumption that
$\Delta_1\cap \Delta_2=\varnothing$, where
$\Delta_2:=[a_1,a_2]$. In the definition of the measure~$s_2$ we follow the standard
convention to the effect that
$d\<\sigma,\sigma_2,\sigma\>:=d\<\sigma,\<\sigma_2,\sigma\>\>$ (for more details,
see~\cite{GoRaSo97}, and
also \cite{ApLy10}, \cite{FiLoLoSo10},~\cite{FiLo11}).
According to what has been said, the three functions $\sigma(z),s_1(z)$ and $s_2(z)$ form
a~(classical) Nikishin system.
This being so, in view of \eqref{10}--\eqref{12}, it is also natural to regard the system of functions $f,f^2,f^3$
as a~Nikishin system, because this system is generated by a~linear combination of three functions,
$\myh{\sigma},\myh{s}_1$ and $\myh{s}_2$, which forms a~Nikishin system.

\section{Proof of Proposition \ref{pro1}}\label{s3}

\subsection{}\label{s3s1}
Given $x_1\in\Delta_1^\circ$, we let $f_j^{+}(x_1)$ and $f_j^{-}(x_1)$, $j=1,2$,
denote the limiting values of the function
$f_j(z)$ as $z=x_1+i\eps\to x_1\in\Delta_1^\circ$, $\eps\to0$, assuming that~$z$
lies, respectively, in the upper half-plane ($\eps>0$) and in the lower half-pane ($\eps<0$).
It is easily seen that
\begin{equation}
f_j^{+}(x_1)=\(A_j-(x_1-i\sqrt{1-x^2_1})\)^{-1/2},\quad
f_j^{-}(x_1)=\(A_j-(x_1+i\sqrt{1-x^2_1})\)^{-1/2}.
\label{33}
\end{equation}
A direct consequence of \eqref{33} is that, for $x_1\in\Delta_1^\circ$,
\begin{align}
\Delta f_j(x_1):
&=(f_j^{+}-f_j^{-})(x_1)\notag\\
&=\frac{\(A_j-(x_1+i\sqrt{1-x^2_1})\)^{1/2}-\(A_j-(x_1-i\sqrt{1-x^2_1})\)^{1/2}}
{\[\(A_j-(x_1+i\sqrt{1-x^2_1})\)\(A_j-(x_1-i\sqrt{1-x^2_1})\)\]^{1/2}}
\notag\\
&=-\frac{2i\sqrt{1-x^2_1}}
{\sqrt{(A_j-x_1)^2+(1-x^2_1)}h_j(x_1)}
\notag\\
&=-\frac{2i\sqrt{1-x^2_1}}
{\sqrt{2A_j(a_j-x_1)}h_j(x_1)},
\label{34}
\end{align}
where
\begin{equation}
h_j(x_1):=\(A_j-(x_1+i\sqrt{1-x^2_1})\)^{1/2}+\(A_j-(x_1-i\sqrt{1-x^2_1})\)^{1/2}
\label{36}
\end{equation}
for $x_1\in\Delta_1^\circ$, $j=1,2$.
Moreover, we have
\begin{align}
(f_j^{+}+f_j^{-})(x_1)
&=\(A_j-(x_1+i\sqrt{1-x^2_1})\)^{-1/2}+\(A_j-(x_1-i\sqrt{1-x^2_1})\)^{-1/2}
\notag\\
&=\frac{h_j(x_1)}{\sqrt{2A_j(a_j-x_1)}}.
\label{35}
\end{align}
It is easily checked that each function $h_j$, which is holomorphic on the interval $\Delta_1^\circ$,
extends holomorphically from this interval to some neighborhood of~$\Delta_1$.
Moreover, $h_j(x_1)\neq0$ for $x_1\in\Delta_1$, and therefore, for $x_1$ from some neighborhood of~$\Delta_1$.
It is also worth noting that the function $f_j^{+}+f_j^{-}$, which is holomorphic on the interval
$\Delta_1^0$, extends holomorphically to some neighborhood of the interval~$\Delta_1$.

We have $f=f_1f_2$, and hence, for $x_1\in\Delta^\circ_1$, using the identity
$$
2\Delta f(x_1):=2(f^{+}-f^{-})(x_1)
=\Delta f_1(x_1)(f_2^{+}+f_2^{-})(x_1)+\Delta f_2(x_1)(f_1^{+}+f_1^{-})(x_1)
$$
and employing relations \eqref{34} and~\eqref{35}, we get
\begin{equation}
2\Delta f(x_1)=
-\frac{i\sqrt{1-x^2_1}}{\sqrt{A_1A_2}\sqrt{(a_1-x_1)(a_2-x_1)}}
\[\frac{h_2(x_1)}{h_1(x_1)}+\frac{h_1(x_1)}{h_2(x_1)}\], \quad x_1\in \Delta_1.
\label{38}
\end{equation}
Moreover, it is also immediate that
$$
\frac{h_1(x_1)}{h_2(x_1)}+
\frac{h_2(x_1)}{h_1(x_1)}>0\quad\text{for}\quad x_1\in\Delta_1.
$$

We have $f(\infty)=1/\sqrt{A_1A_2}$ by definition~\eqref{8} of the function~$f$, and hence,
applying Cauchy's theorem to the function~$f$, we get the following representation
\begin{equation}
f(z)-\frac1{\sqrt{A_1A_2}}
=\frac1{2\pi i}\int_{\gamma_1}\frac{f(t)}{t-z}\,dt,
\quad z\in\ext\gamma_1,
\label{38.2}
\end{equation}
where $\gamma_1$ is an arbitrary closed Jordan curve separating the interval
$\Delta_1$ from the infinity point and containing the point~$z$ in the unbounded component
$\ext\gamma_1$ of its complement $\myo\CC\setminus\gamma_1$; we assume that the curve~$\gamma_1$
has positive orientation relative to $\ext\gamma_1$. From~\eqref{38.2} it easily follows that
\begin{equation}
f(z)-\frac1{\sqrt{A_1A_2}}
=\frac1{2\pi i}\int_{\Delta_1}\frac{\Delta f(x_1)}{x_1-z}\,dx
=\int_{\Delta_1}\frac{d\sigma(x_1)}{z-x_1}=\myh{\sigma}(z),
\label{40}
\end{equation}
where, for $x_1\in\Delta_1$,
\begin{align}
d\sigma(x_1)
&=-\frac1{2\pi i}\Delta f(x_1)\,dx
\notag\\
&=
\frac{\sqrt{1-x^2_1}}{4\pi\sqrt{A_1A_2}\sqrt{(a_1-x_1)(a_2-x_1)}}
\[\frac{h_2(x_1)}{h_1(x_1)}+\frac{h_1(x_1)}{h_2(x_1)}\]\,dx_1.
\label{41}
\end{align}
Using~\eqref{40} and~\eqref{41}, this establishes
\begin{equation}
f(z)=\frac1{\sqrt{A_1A_2}}+\myh{\sigma}(z),
\quad z\in D_1,
\notag
\end{equation}
thereby proving representations \eqref{10} and~\eqref{13}.

\subsection{}\label{s3s2}
We set $\rho_1(x_1):=f^{+}(x_1)+f^{-}(x_1)=(f_1^{+}f_2^{+}+f_1^{-}f_2^{-})(x_1)$,
$x_1\in\Delta_1^\circ$. It is easily seen (see~\eqref{8} and~\eqref{35}) that the
function $\rho_1\in\HH(\Delta^\circ_1)$ extends holomorphically from the interval $\Delta_1^\circ$ to some
neighborhood of the interval~$\Delta_1$. Moreover, the function $\rho_1$ is holomorphic on the domain
$D_2:=\myo\CC\setminus\Delta_2$ and can be represented in this domain as
\begin{align}
\rho_1(z)
&=\[
\(A_1-(z-(z^2-1)^{1/2})\)\(A_2-(z-(z^2-1)^{1/2})\)\]^{-1/2}
\notag\\
&\quad+\[\(A_1-(z+(z^2-1)^{1/2})\)\(A_2-(z+(z^2-1)^{1/2})\)\]^{-1/2}
\notag\\
&=\[\(A_1-\frac1{\pfi(z)}\)\(A_2-\frac1{\pfi(z)}\)\]^{-1/2}+
\[\bigl(A_1-\pfi(z)\bigr)\bigl(A_2-\pfi(z)\bigr)\]^{-1/2}.
\label{43}
\end{align}
Given $x_2\in\Delta_2^\circ$, we denote by $\rho_1^{+}(x_2)$ the limiting values of the function
$\rho_1(z)$ as $z\to x_2$ assuming that $z$~lies in the upper half-plane, and denote by
$\rho_1^{-}(x_2)$ the limiting values of $\rho_1(z)$ as $z\to x_2$ assuming that $z$~lies in the lower half-plane.
Using~\eqref{43},
\begin{equation}
\Delta \rho_1(x_2):=(\rho_1^{+}-\rho_1^{-})(x_2)
=\frac{-2i}{\sqrt{(\pfi(x_2)-A_1)(A_2-\pfi(x_2))}}, \quad x_2\in\Delta_2^\circ.
\label{44}
\end{equation}
We have $\rho_1(\infty)=1/\sqrt{A_1A_2}$. Hence, by~\eqref{43}
\begin{equation}
\rho_1(z)-\frac1{\sqrt{A_1A_2}}=\frac1{2\pi i}\int_{\gamma_2}\frac{\rho_1(t)\,dt}{t-z}
=\frac1{2\pi i}\int_{\Delta_2}\frac{\Delta \rho_1(x_2)\,dx_2}{x_2-z},
\label{45}
\end{equation}
where $\gamma_2$ is an arbitrary negatively oriented closed Jordan curve
separating the interval~$\Delta_2$ from the
infinity point; the point~$z$ lies in that connected component of
$\myo\CC\setminus\gamma_2$ which contains the infinity point.

From~\eqref{44} and~\eqref{45} we see that
\begin{equation}
\rho_1(z)=\frac1{\sqrt{A_1A_2}}+\myh{\sigma}_2(z),\quad z\in D_2,
\label{46}
\end{equation}
where
\begin{equation}
d\sigma_2(x_2):=-\frac1{2\pi i}\Delta \rho_1(x_2)\,dy
=\frac1{\pi}\frac{dx_2}{\sqrt{(\pfi(x_2)-A_1)(A_2-\pfi(x_2)}},\quad x_2\in\Delta_2^\circ.
\label{47}
\end{equation}
So, we have
$$
\rho_1(z):=(f^{+}+f^{-})(z)=\frac1{\sqrt{A_1A_2}}+\myh{\sigma}_2(z),
$$
where $\sigma_2$ is the positive measure with support in~$\Delta_2$ defined by
representation~\eqref{47}. Hence, for $x_1\in\Delta_1$,
\begin{equation}
\frac{\Delta f^2}{\Delta f}(x_1)=(f^{+}+f^{-})(x_1)
=\frac1{\sqrt{A_1A_2}}+\myh{\sigma}_2(x_1).
\label{47.2}
\end{equation}
As a result (see~\eqref{41}), we have, for $x_1\in\Delta_1$,
\begin{align}
\Delta f^2(x_1)\,dx_1
&=\(\frac1{\sqrt{A_1A_2}}+\myh{\sigma}_2(x_1)\)\Delta f(x_1)\,dx_1
\notag\\
&=-\(\frac1{\sqrt{A_1A_2}}+\myh{\sigma}_2(x_1)\)2\pi i\,d\sigma(x_1).
\label{48}
\end{align}
Since $f^2(\infty)=1/(A_1A_2)$, it follows from~\eqref{48} that
\begin{align}
f^2(z)-\frac1{A_1A_2}&=
\frac1{2\pi i}\int_{\gamma_1}\frac{f^2(t)}{t-z}\,dt
=-\frac1{2\pi i}\int_{\Delta_1}\frac{\Delta f^2(x_1)}{z-x_1}\,dx_1
\notag\\
&=\frac1{\sqrt{A_1A_2}}\myh{\sigma}(z)
+\int_{\Delta_1}\frac{\myh{\sigma}_2(x_1)\,d\sigma(x_1)}{z-x_1}
=\frac1{\sqrt{A_1A_2}}\myh{\sigma}(z)+\myh{s}_1(z),
\notag
\end{align}
where $s_1=\<\sigma,\sigma_2\>$, $\supp s_1=\Delta_1$. Therefore,
\begin{equation}
f^2(z)=\frac1{A_1A_2}+\frac1{\sqrt{A_1A_2}}\myh{\sigma}(z)-\myh{s}_1(z),
\quad z\in D_1.
\notag
\end{equation}
This completes the proof of representations \eqref{11} and~\eqref{14}.

\subsection{}\label{s3s3}
We now set
\begin{equation}
\rho_2(x_1):=\frac{\Delta f^3(x_1)}{\Delta f(x_1)}, \quad x_1\in \Delta_1^\circ.
\label{51}
\end{equation}
Given $x_1\in\Delta_1^\circ$, we have
\begin{align}
f^{+}(x_1)&=\[\(A_1-(x_1-i\sqrt{1-x^2_1})\)\(A_2-(x_1-i\sqrt{1-x^2_1})\)\]^{-1/2},
\notag\\
f^{-}(x_1)&=\[\(A_1-(x_1+i\sqrt{1-x^2_1})\)\(A_2-(x_1+i\sqrt{1-x^2_1})\)\]^{-1/2},
\notag
\end{align}
and hence,
\begin{align}
(f^{+})^2(x_1)&=
\frac1{\(A_1-(x_1-i\sqrt{1-x^2_1})\)\(A_2-(x_1-i\sqrt{1-x^2_1})\)},
\notag\\
(f^{-})^2(x_1)&
=\frac1{\(A_1-(x_1+i\sqrt{1-x^2_1})\)\(A_2-(x_1+i\sqrt{1-x^2_1})\)}.
\notag
\end{align}
It follows that the functions $f^{+}f^{-}$ and $(f^{+})^2+(f^{-})^2$
extend analytically from the interval $\Delta_1^\circ$ to the domain $D_2$. Consequently,
the function $\rho_2$, which is given by representation \eqref{51}, extends holomorphically to the domain~$D_2$.
Moreover,
\begin{gather}
f^{+}(x_1)f^{-}(x_1)=\frac1{\sqrt{2A_1(a_1-x_1)}\sqrt{2A_2(a_2-x_1)}}
=\frac1{2\sqrt{A_1(a_1-x_1)A_2(a_2-x_1)}},
\notag\\
\begin{aligned}
(f^{+})^2(x_1)+(f^{-})^2(x_1)
&=
\frac{\(A_1-(x_1+i\sqrt{1-x^2_1})\)\(A_2-(x_1+i\sqrt{1-x^2_1})\)}
{4A_1(a_1-x_1)A_2(a_2-x_1)}
\notag\\
&\quad+\frac{\(A_1-(x_1-i\sqrt{1-x^2_1})\)\(A_2-(x_1-i\sqrt{1-x^2_1})\)}
{4A_1(a_1-x_1)A_2(a_2-x_1)}
\end{aligned}
\notag
\end{gather}
for $x_1\in\Delta_1^\circ$. Since $1<a_1<a_2$, it can be easily shown that
$((f^{+})^2+(f^{-})^2+f^{+}f^{-})(x_1)>0$ for
$x_1\in\Delta_1$. So, using the definition of the function $\pfi(z)$ and
employing the identity
$$
\frac{\Delta f^3(x_1)}{\Delta f(x_1)}=((f^{+})^2+(f^{-})^2)(x_1)
+(f^{+}f^{-})(x_1),
$$
where $x\in\Delta_1^\circ$, we arrive at the explicit representation
\begin{align}
\rho_2(z)&=
\frac1{(A_1-1/\pfi(z))(A_2-1/\pfi(z))}+
\frac1{(A_1-\pfi(z))(A_2-\pfi(z))}
\notag\\
&\quad+
\[(A_1-1/\pfi(z))(A_2-1/\pfi(z))(A_1-\pfi(z))(A_2-\pfi(z))\]^{-1/2}
\label{55}
\end{align}
for the function $\rho_2\in\HH(D_2)$, where $z\in D_2$.
Moreover, $\rho_2(x_1)>0$ for $x_1\in\Delta_1$, $\rho_2(\infty)=1/(A_1A_2)$, and for
$x_2\in\Delta_2^\circ$, we have
\begin{align}
\Delta \rho_2(x_2):&=\rho_2^{+}(x_2)-\rho_2^{-}(x_2)
\notag\\
&=\frac{-2i}{\sqrt{(A_1-1/\pfi(x_2))(A_2-1/\pfi(x_2))
(\pfi(x_2)-A_1)(A_2-\pfi(x_2))}}.
\label{56}
\end{align}
Therefore,
\begin{align}
\rho_2(z)
&=\rho_2(\infty)
+\frac1{2\pi i}\int_{\gamma_2}\frac{\rho_2(t)\,dt}{t-z}
\notag\\
&=\frac1{A_1A_2}
+\frac1{2\pi i}\int_{a_1}^{a_2}\frac{\Delta \rho_2(x_2)\,dx_2}{x_2-z},
\label{57}
\end{align}
where $z\in D_2$, $\gamma_2$ is an arbitrary closed Jordan curve
separating the interval $\Delta_2$ from the point~$z$ and from the
infinity point and which is positively oriented with respect to the
domain containing the point~$z$. The following representation for the
function $\rho_2(z)$ is a~direct consequence of \eqref{56}
and~\eqref{57}. We have
\begin{equation}
\rho_2(z)=\frac1{A_1A_2}+\myh{\sigma}_3(z),
\label{58}
\end{equation}
where $\sigma_3$ is a positive measure supported on the interval $\Delta_2$,
$\supp\sigma_3=\Delta_2$, and moreover,
\begin{equation}
d\sigma_3(x_2)=
\frac1\pi\frac{dx_2}
{\sqrt{(A_1-1/\pfi(x_2))(A_2-1/\pfi(x_2))(\pfi(x_2)-A_1)(A_2-\pfi(x_2))}},
\label{59}
\end{equation}
$x_2\in\Delta_2^\circ$.

So, for $z\in D_2$ we have
\begin{equation}
\frac{\Delta f^3}{\Delta f}(z)=\rho_2(z)=\frac1{A_1A_2}+\myh{\sigma}_3(z).
\label{60}
\end{equation}
Hence, in view of \eqref{47} it follows from~\eqref{59} that
\begin{equation}
d\sigma_3(x_2)=\rho_3(x_2)\,d\sigma_2(x_2),\quad x_2\in\Delta_2^\circ,
\label{61}
\end{equation}
where
$$
\rho_3(x_2):=\frac1{\sqrt{(A_1-1/\pfi(x_2))(A_2-1/\pfi(x_2))}},
\quad x_2\in\Delta_2^\circ.
$$
The function $\rho_3$ extends holomorphically from the interval $\Delta_2^\circ$ to the domain
$D_1$. Fur\-ther\-more, it is clear that $\rho_3(z)\equiv f(z)$, $z\in D_1$.
So, by~\eqref{10}
\begin{equation}
\rho_3(z)=\frac1{\sqrt{A_1A_2}}+\myh{\sigma}(z),\quad z\in D_1,
\label{62}
\end{equation}
where the measure $\sigma$ is given by representation~\eqref{41}.

From~\eqref{60},~\eqref{61} and~\eqref{62} it follows that, for $x_1\in\Delta_1^\circ$,
\begin{align}
\frac{\Delta f^3(x_1)}{\Delta f(x_1)}
&=\frac1{A_1A_2}+\int_{a_1}^{a_2}\frac{\rho_3(x_2)\,d\sigma_2(x_2)}{x_1-x_2}
\notag\\
&=\frac1{A_1A_2}+\frac1{\sqrt{A_1A_2}}\myh{\sigma}_2(x_1)
+\int_{a_1}^{a_2}\frac{\myh{\sigma}(x_2)\,d\sigma_2(x_2)}{x_1-x_2}
\notag\\
&=\frac1{A_1A_2}+\frac1{\sqrt{A_1A_2}}\myh{\sigma}_2(x_1)+\myh{s}(x_1),
\label{63}
\end{align}
where the measure $s$ is defined as $s:=\<\sigma_2,\sigma\>$,
$\supp{s}=\supp{\sigma_2}=\Delta_2$.
We have $f^3(\infty)=1/\sqrt{(A_1A_2)^3}$, and hence, by Cauchy's formula,
\begin{equation}
f^3(z)-\frac1{\sqrt{(A_1A_2)^3}}
=\frac1{2\pi i}\int_{\gamma_1}\frac{f^3(t)\,dt}{t-z}
=\frac1{2\pi i}\int_{-1}^1\frac{\Delta f^3(x_1)\,dx_1}{x_1-z}.
\label{64}
\end{equation}
So, using \eqref{60},~\eqref{63} and~\eqref{64},
\begin{align}
f^3(z)-\frac1{\sqrt{A_1A_2}}
&=\frac1{A_1A_2}\cdot\frac1{2\pi i}\int_{-1}^1\frac{\Delta f(x_1)\,dx_1}{x_1-z}
+\frac1{\sqrt{A_1A_2}}\cdot\frac1{2\pi
i}\int_{-1}^1\frac{\myh{\sigma}_2(x_1)\Delta f(x_1)\,dx_1}{x_1-z}
\notag\\
&\quad-\frac1{2\pi i}\int_{-1}^1\frac{\myh{s}(x_1)\Delta f(x_1)\,dx_1}{x_1-z}
\notag\\
&=\frac1{A_1A_2}\myh{\sigma}(z)
+\frac1{\sqrt{A_1A_2}}\int_{-1}^1\frac{\myh{\sigma}_2(x_1)\,d\sigma(x_1)}{z-x_1}
+\int_{-1}^1\frac{\myh{s}(x_1)\,d\sigma(x_1)}{z-x_1}.
\label{65}
\end{align}
Finally, from~\eqref{65} and the definition of the measure $s=\<\sigma_2,\sigma\>$ we have the
representation
$$
f^3(z)=\frac1{\sqrt{(A_1A_2)^3}}
+\frac1{A_1A_2}\myh{\sigma}(z)+\frac1{\sqrt{A_1A_2}}\myh{s}_1(z)
+\myh{s}_2(z),
$$
where $s_1=\<\sigma,\sigma_2\>$, $s_2=\<\sigma,\sigma_2,\sigma\>$,
$\supp{s_j}=[-1,1]$, $j=1,2$.

This proves \eqref{12}, and therefore, Proposition~\ref{pro1}.

\begin{remark}\label{rem2}
The relations
\begin{equation}
f(z)-\frac1{\sqrt{A_1A_2}}=\myh{\sigma}(z),\quad
f^2(z)-\frac1{\sqrt{A_1A_2}}f(z)=\myh{s}_1(z),\quad
f^3(z)-\frac1{\sqrt{A_1A_2}}f^2(z)=\myh{s}_2(z)
\label{14.2}
\end{equation}
are immediate consequences of \eqref{10}--\eqref{12}.
\end{remark}

\begin{conj}\label{con1}
Let $f$ be a~function from the class $\sZ$ of the form
\begin{equation}
f(z)=\(\frac{A_1-1/\pfi(z)}{A_2-1/\pfi(z)}\)^\alpha,
\label{71}
\end{equation}
where $1<A_1<A_2$ and $\alpha\in\RR\setminus\QQ$. Then $f$ is a~Markov function
and, for any $n\in\NN$, the system $f,f^2,\dots,f^n$ is a~Nikishin system.
\end{conj}

\begin{remark}\label{rem3}
In accordance with representation~\eqref{8} all branch points of the function~$f$ are
of second order, and hence in view of the above Proposition~\ref{pro1}
the results of \cite{ApBoYa17} and~\cite{LoVa18} on the asymptotics of Hermite--Pad\'e polynomials
apply to the system of functions $f,f^2$, and the results of~\cite{KoPaSuCh17}, to the system of
functions $f,f^2,f^3$. It is very likely
that by appropriately transforming the independent variable
(see, for example, \cite[\S~5]{Sue18c}), which was carried out in representation~\eqref{8},
and multiplying some resulting functions it might be possible
to obtain those exotic, as they may seem, Nikishin systems on star-like sets which have been considered in
\cite{BaGeLo18} and~\cite{LoMi18}. In other words, there is a~hope that examples of Nikishin systems
of such kind can be found in the form $f,f^2,\dots,f^n$.
\end{remark}


\end{document}